
\baselineskip=16pt
\nopagenumbers

\def\RRR{{\bf R}}
\def\IPxyl{P^{x,y}_t}
\def\leb{\mathop{\rm Leb}}
\def\nin{\noindent}
\def\ref#1{{\rm [}{\bf #1}{\rm ]}}   
\def\nref#1#2{{\rm [}{\bf #1}{\rm ;\ #2]}}
\outer\def\proclaim#1{\medbreak\noindent\bf\ignorespaces
   #1\unskip.\enspace\sl\ignorespaces}
\outer\def\endproclaim{\par\ifdim\lastskip<\medskipamount\removelastskip
   \penalty 55 \fi\medskip\rm}
\def\comp{\leavevmode
     \raise.2ex\hbox{${\scriptstyle\mathchar"020E}$}}
\def\prf{\nin{\it Proof. }}
\def\rect#1#2#3{\raise .1ex\vbox{\hrule height.#3pt
   \hbox{\vrule width.#3pt height#2pt \kern#1pt\vrule width.#3pt}
        \hrule height.#3pt}}
\def\qed{$\hskip 5pt\rect364$} 

\def\FF{{\cal F}}\def\EE{{\cal E}}\def\BB{{\cal B}}\def\II{{\cal I}}

\def\BL{BL97} 
\def\BIL{B68} 
\def\DYN{D80}
\def\FPY{FPY93} 
\def\FOL{F90} 
\def\GET{G75}
\def\GSexc{GS82}
\def\GSnswd{GS84}
\def\KW{KW65}
\def\KUN{K76}
\def\SHA{S88}
\def\WIT{W86}
\def\YAN{Y88}

\topinsert
\vskip1truein
\endinsert

\centerline{\bf Markov Processes with Identical Bridges}
\bigskip 
\centerline{by}
\medskip
\centerline{P. J. Fitzsimmons}

\centerline{Department of Mathematics}
\centerline{University of California, San Diego}
\centerline{9500 Gilman Drive}
\centerline{La Jolla, CA  92093-0112\ USA}
\centerline{\tt pfitzsim@ucsd.edu}
\bigskip\bigskip\bigskip

\centerline{\bf ABSTRACT}
Let $X$ and $Y$ be time-homogeneous Markov processes with  common state space
$E$, and assume that the transition kernels of $X$ and $Y$ admit densities with
respect to suitable reference measures. We show that if there is a time $t>0$
such that, for each $x\in E$, the conditional distribution of $(X_s)_{0\le s\le
t}$, given $X_0=x=X_t$, coincides with the conditional distribution of $(Y_s)_{0\le s\le
t}$, given $Y_0=x=Y_t$, then the infinitesimal generators of $X$ and $Y$ are
related by $L^Yf=\psi^{-1}L^X(\psi f)-\lambda f$, where $\psi$ is an eigenfunction
of $L^X$ with eigenvalue $\lambda\in{\bf R}$. Under an additional continuity
hypothesis, the same conclusion obtains assuming merely that $X$ and $Y$ share a
``bridge'' law for one triple $(x,t,y)$. Our work entends and clarifies a recent
result of I.\ Benjamini and S.\ Lee.
\medskip

\bigskip\bigskip\bigskip

\nin{\it Running Title:\/} Equal Bridges
\medskip

\nin{\it Key words and phrases.\/} Bridge law, eigenfunction, transition density.
\medskip

\nin{\it 1990 AMS Subject classification.\/} Primary: 60J25; secondary
60J35.
 
\eject
\baselineskip=16pt
\pageno=1
\footline{\hfill{\rm\folio}\hfill}

\nin{\bf 1. Introduction}
\medskip

\nin Let $X=(X_t,P^x)$ and $Y=(Y_t,Q^x)$ be  non-explosive
regular Markov diffusion processes in ${\bf R}$. Let $P^{x,y}_t$ denote the
conditional law of
$(X_s)_{0\le s\le t}$ given $X_0=x$, $X_t=y$. Let $Q^{x,y}_t$ denote the
analogous ``bridge'' law for $Y$. Recently, Benjamini \& Lee \ref{\BL} 
proved the following result.

\proclaim{(1.1) Theorem} Suppose that $X$ is standard Brownian motion and that
$Y$ is a weak solution of the stochastic differential equation
$$
dY_t=dB_t+\mu(Y_t)\,dt,
\leqno(1.2)
$$
where $B$ is standard Brownian motion and the drift $\mu$ is bounded and twice
continuously differentiable. If $Q^{x,x}_t=P^{x,x}_t$ for all
$x\in{\bf R}$ and all
$t>0$, then either (i) $\mu(x)\equiv k$ or (ii)
$\mu(x)=k\tanh(kx+c)$, for some real constants $k$ and $c$.
\endproclaim

Our aim in this paper is to generalize this theorem in two ways.

Firstly, we allow $X$ and $Y$ to be general strong Markov processes with
values in an abstract state space $E$.  We  require that $X$ and $Y$ have
dual processes with respect to suitable reference measures, and that $X$ and $Y$
admit  transition densities with respect to these reference measures.
(These conditions are met by all regular 1-dimensional diffusions without
absorbing boundary points.)

Secondly, under an additional continuity condition, we  show that the equality of
$Q^{x,y}_t$ and $P^{x,y}_t$ for a {\it single\/} choice of the triple $(x,t,y)$ is enough
to imply that
$Q^{x,y}_t=P^{x,y}_t$ for {\it all\/} $(x,t,y)\in E\times ]0,\infty[\times  E$.
We provide a simple example illustrating what can go wrong when the continuity
condition fails to hold.

The conclusion of Theorem (1.1) is more transparently  stated as follows. Given a
drift  $\mu$  define 
$\psi(x):=\exp\int_0^x \mu(y)\,dy$. Then  $\mu$ satisfies
 the conclusion of Theorem (1.1) if and only if
$$
{\textstyle{1\over 2}}\psi''(x)=\lambda\,\psi(x),\qquad \forall\,  x\in{\bf R},
$$
where $\lambda:=k^2/2$. Thus, Theorem (1.1) can be stated as follows: If $X$ is
Brownian motion and if $Y$ is ``Brownian motion with drift $\mu$,'' then
$X$ and
$Y$  have common bridge laws if and only if
$\mu$ is the logarithmic derivative of a strictly positive  eigenfunction of the local
infinitesimal generator of $X$, in which case the  laws of $X$ and $Y$ are
related by
$$
{dQ^x\over dP^x}\Big|_{\FF_t}=e^{-\lambda t}{\psi(X_t)\over\psi(X_0)}.
\leqno(1.3)
$$

Theorem (1.1) and our extensions of it depend crucially on the existence of a
``reference'' measure dominating the transition probabilities of $X$ and $Y$. This fact
is amply demonstrated by the work of H.\ F\"ollmer in
\ref{\FOL}.   
Let $E$ be the Banach space of continuous maps of $[0,1]$ into
${\bf R}$ that vanish at 0, and let $m$ denote Wiener measure on the Borel subsets of
$E$. Let $X=(X_t,P^x)$ be the associated Brownian motion in $E$; that is, the $E$-valued
diffusion with transition semigroup given by
$$
P_t(x,f):=\int_E f(x+\sqrt{t}y)\,m(dy).
$$
This semigroup admits no reference measure; indeed $P_t(x,\cdot\,)\perp P_t(y,\cdot\,)$
unless $x-y$ is an element of the Cameron-Martin space $H$, consisting of those elements
of $E$ that are absolutely continuous and possess a square-integrable derivative.
Now given $z\in E$, let $Y=(Y_t,Q^x)$ be Brownian motion in $E$ with drift
$z$. By this we mean the $E$-valued diffusion with transition semigroup
$$
Q_t(x,f):=\int_E f(x+tz+\sqrt{t}y)\,m(dy).
$$
Given
$(x,t,y)\in E\times]0,\infty[\times E$,  let
$P^{x,y}_t$ be the $P^0$-distribution of the process $\{x+X_s+(s/t)(y-x-X_t): 0\le s\le
t\}$. Evidently, (i) $(x,y)\mapsto P^{x,y}_t$ is weakly continuous, (ii)
$P^{x,y}_t(X_t=y)=1$, and (iii) 
$\{P^{x,y}_t:y\in E\}$ is a regular version of the family of conditional
distributions
$Q^x(\{X_s;0\le s\le t\}\in\cdot\,|\,X_t=y)$, regardless of the choice of $z\in E$.
In other words, $X$ and
$Y$ have common bridge laws. However, the laws of $X$ and $Y$ are mutually absolutely
continuous  (as in (1.3)) if and only if $z\in H$.

Before stating our results we
describe the context in which we shall be working. Let
$X=(X_t,P^x)$ now denote a strong Markov process with cadlag paths and infinite lifetime.
We assume that the state space
$E$ is homeomorphic to a Borel subset of some compact metric space, and that the
transition semigroup
$(P_t)_{t\ge 0}$ of
$X$ preserves Borel measurability and is without branch points. In other words, $X$ is a
Borel right processes with cadlag paths and infinite lifetime; see \ref{\GET, \SHA}. The
process
$X$ is realized as the coordinate process
$X_t:\omega\mapsto\omega(t)$ on the sample space
$\Omega$ of all cadlag paths from $[0,\infty[$ to $E$. The probability
measure $P^x$ is the law of
$X$ under the initial condition $X_0=x$.  We write $(\FF_t)_{t\ge 0}$ for the
natural (uncompleted) filtration of $(X_t)_{t\ge 0}$ and $(\theta_t)_{t\ge 0}$
for the  shift operators on $\Omega$: $X_s\comp\theta_t=X_{s+t}$. 

In addition, we assume the existence of transition densities with respect to a
reference measure and (for technical reasons) the existence of a dual process.
(The duality hyothesis (1.4) can be replaced by conditions ensuring the existence of a
nice Martin exit boundary for the space-time process $(X_t,r+t)_{t\ge 0}$; see
\ref{\KW}.) 

Let $\EE$ denote
the Borel $\sigma$-algebra on $E$.
\medskip

\nin{\bf (1.4) Hypothesis.} (Duality) There is a $\sigma$-finite
measure
$m^X$ on $(E,\EE)$ and   a second $E$-valued Borel right Markov
process
$\hat X$, with cadlag paths and infinite lifetime, such that the
semigroup $(\hat P_t)$ of
$\hat X$ is in duality with $(P_t)$ relative to $m^X$:
$$
\int _E f(x) P_tg(x)\,m^X(dx)=\int_E \hat P_tf(x) g(x)\,m^X(dx),
\leqno(1.5)
$$
for all $t>0$ and all positive $\EE$-measurable functions $f$ and $g$.
\medskip

\nin{\bf (1.6) Hypothesis.} (Transition densities)  There is an
$\EE\otimes\BB_{]0,\infty[}\otimes \EE$-measurable function
$(x,t,y)\mapsto p_t(x,y)\in]0,\infty[$ such that
$$
P^x(f(X_t))=P_tf(x)=\int_E p_t(x,y)\,f(y)\,m^X(dy),\qquad  \forall\,t> 0,
\leqno(1.7)
$$
and
$$
\hat P^x(f(X_t))=\hat P_tf(x)=\int_E p_t(y,x)\,f(y)\,m^X(dy),\qquad  \forall\,t> 0,
\leqno(1.8)
$$
for any bounded $\EE$-measurable function $f$. Furthermore, we assume that the
Chapman-Kolmogorov identity holds:
$$
p_{t+s}(x,y)=\int_E p_t(x,z)p_s(z,y)\,m^X(dz),\qquad\forall\,s,t>0, x,y\in E.
\leqno(1.9)
$$
\medskip

Hypothesis (1.6) implies that $m^X(U)>0$ for every non-empty finely open
subset of $E$.

When (1.4) is in force, the existence and uniqueness of a (jointly measurable) transition
density function
$p_t(x,y)$ such that (1.7)--(1.9) hold is guaranteed by the apparently weaker condition:
$P_t(x,\cdot)\ll m^X$, $\hat P_t(x,\cdot)\ll m^X$ for all $x\in E$, $t>0$. 
See, for example,  \ref{\DYN, \WIT, \YAN}. For more discussion of processes with ``dual
transition densities,'' see \nref{\GSexc}{\S 3}.  

Let $Y=(Y_t,Q^x)$ be a second $E$-valued Borel right Markov process with cadlag paths and
infinite lifetime. The process $Y$ is assumed to satisfy all of the conditions imposed
on $X$ above. In particular, we can (and do) assume that $Y$ is realized as the
coordinate process on $\Omega$. The transition semigroup of
$Y$ is denoted
$(Q_t)_{t\ge 0}$ and we use
$m^Y$ and $q_t(x,y)$ to denote the reference measure and transition density function for
$Y$. (The bridge laws $P^{x,y}_t$
and $Q^{x,y}_t$ for $X$ and $Y$ will be discussed in more detail in section 2.)
 
In what follows, the prefix ``co-'' refers to the dual process $\hat X$ (or $\hat Y$). 

\proclaim{(1.10) Theorem} Let $X$ and $Y$ be  strong Markov processes as
described above, satisfying Hypotheses (1.4) and (1.6). Suppose  there exists
 $t_0>0$ such that 
$Q^{x,x}_{t_0}=P^{x,x}_{t_0}$ for all $x\in E$. Then

\itemitem{(a)} $P^x|_{\FF_t}\sim Q^x|_{\FF_t}$ and $\hat P^y|_{\FF_t}\sim
\hat Q^y|_{\FF_t}$,  for all $x\in E$, $y\in E$, and $t>0$;

\itemitem{(b)} There exist a constant
$\lambda\in\RRR$, a Borel finely continuous function
$\psi:E\to]0,\infty[$, and a  Borel co-finely continuous function
$\hat\psi:E\to]0,\infty[$ such that for all $t>0$,
\vskip-15pt
$$
P_t\psi(x)=e^{\lambda t}\,\psi(x),\qquad\forall\,x\in E,
\leqno(1.11)
$$ 
$$
\hat P_t\hat\psi(x)=e^{\lambda
t}\,\hat\psi(x),\qquad\forall\,x\in E,
\leqno(1.12)
$$ 
$$
Q^x|_{\FF_t}=e^{-\lambda t}{\psi(X_t)\over\psi(X_0)}\,P^x|_{\FF_t},\qquad\forall\,x\in E,
\leqno(1.13)
$$
$$
\hat Q^x|_{\FF_t}=e^{-\lambda t}{\hat\psi(X_t)\over\hat\psi(X_0)}\,\hat P^x|_{\FF_t},
\qquad\forall\,x\in E.
\leqno(1.14)
$$
The function $\psi\hat\psi$ is a Borel version
of the Radon-Nikodym derivative $dm^Y/dm^X$.

\itemitem{(c)} $Q^{x,y}_t=P^{x,y}_t$ for {\it all\/} $(x,t,y)\in E\times]0,\infty[\times
E$;
\endproclaim

\nin{\bf (1.15) Remarks.}

(i) Given  
functions
$\psi$ and $\hat \psi$ as in (1.11) and (1.12), the right sides of (1.13) and
(1.14) determine the
laws of  Borel right Markov processes $Y^*$ and $\hat Y^*$ on
$E$. It is easy to check that $Y^*$ and $\hat Y^*$ are in duality with respect to the
measure $\psi\hat\psi\cdot m^X$, that Hypotheses (1.4) and (1.6) are satisfied, and that
$Y^*$ (resp.\ $\hat Y^*$) has the same bridge laws as $X$ (resp.\ $\hat X$).

(ii) As noted  earlier, any one-dimensional regular diffusion without 
absorbing boundaries satisfies Hypotheses (1.4) and (1.6). Such a diffusion is
self-dual with respect to its speed measure, which serves as the  reference measure.
Moreover, the transition density function of such a diffusion is jointly continuous in
$(x,t,y)$. See
\nref{IM}{pp.\ 149--158}.
\medskip

\proclaim{(1.16) Theorem} Let $X$ and $Y$ be right Markov processes as described before
the statement of Theorem (1.10). Suppose, in addition to (1.4) and (1.6), that for each
$t>0$ the transition  density functions
$p_t(x,y)$ and $q_t(x,y)$ are separately continuous in  the spatial variables $x$ and
$y$. If there is a triple $(x_0,t_0,y_0)\in E\times]0,\infty[\times E$ such that
$P^{x_0,y_0}_{t_0}=Q^{x_0,y_0}_{t_0}$, then the conclusions (a), (b), and (c) of
Theorem (1.10) remain true.
\endproclaim

\nin{\bf (1.17) Remark.} 
Let us suppose that
$X$ is a real-valued regular diffusion  on its natural scale, and that its speed measure
$m^X$ admits a strictly positive density $\rho$ with respect to Lebesgue measure. 
Let
$L^X$ denote the local infinitesimal generator of
$X$. Then (1.11) implies  $L^X\psi=\lambda\psi$, or more explicitly
$$
{1\over \rho(x)}\psi''(x)=\lambda\psi(x).
$$
Moreover, (1.13) means that the transition semigroups of $X$ and $Y$ are
related by 
$$
Q_t(x,dy)=\exp(-\lambda t)[\psi(y)/\psi(x)]P_t(x,dy).
$$
From this it follows that the (local)
infinitesimal generators of $X$ and $Y$ are related by
$$
L^Yf(x)=L^Xf(x)+{2\,\mu(x)\over \rho(x)}\cdot f'(x),
\leqno(1.18)
$$
where $\mu:=(\log\psi)'$.
When
$X$ is standard Brownian motion (so that
$\rho(x)\equiv 2$), the right side of (1.18) is the infinitesimal generator of
any weak solution of (1.2). By Remark (1.15)(ii), the additional condition
imposed in Theorem (1.16) is met in the present situation. Consequently, Theorem (1.16)
implies that the conclusion of Theorem (1.1) is true once we know that the
$(x_0,t_0,y_0)$-bridge law of $Y$ is a Brownian bridge, for {\it one\/} triple
$(x_0,t_0,y_0)$
\medskip

Without some sort of additional condition as in Theorem (1.16), there may be an
exceptional set in the conclusions (a)--(c).  Recall that a
Borel set
$N\subset E$ is
$X$-polar if and only if
$P^x(X_t\in N$ for some
$t>0)=0$ for all $x\in E$.
\medskip

\nin{\bf (1.19) Example.} The state space in this example will be the real line 
${\bf R}$. Let $Z=(Z_t,R^x)$ be a 3-dimensional Bessel process, with state space
$[0,\infty[$. (Under $R^x$, $(Z_t)_{t\ge 0}$ has the same law as the radial part of a
standard 3-dimensional Brownian motion started at $(x,0,0)$.) We assume that the
probability space on which $Z$ is realized is rich enough to support an independent
unit-rate Poisson process $(N(t))_{t\ge 0}$. The process $X$ is presented
(non-canonically) as follows:
$$
X_t:=\cases{(-1)^{N(t)}Z_t,&if $X_0\ge 0$;\cr (-1)^{N(t)+1}Z_t, &if $X_0<0$,\cr}
$$
whereas $Y$ is presented as
$$
Y_t:=\cases{(-1)^{N(t)}Z_t,&if $Y_0> 0$;\cr (-1)^{N(t)+1}Z_t, &if $Y_0\le 0$.\cr}
$$
Both $X$ and $Y$ are Borel right Markov processes satisfying (1.4) and (1.6); indeed, both
processes are self-dual with respect to the reference measure $m(dx):=x^2\,dx$. The
singleton $\{0\}$ is a polar set for both processes. If neither $x$ nor $y$ is equal to
0, then  $P^{x,y}_t=Q^{x,y}_t$ for all $t>0$. However, $P^{0,y}_t$ and $Q^{0,y}_t$ are
different for all $y\in{\bf R}$ and $t>0$, because
$$
P^{0,y}_t(X_s>0\hbox{ for all small }s)=
Q^{0,y}_t(X_s<0\hbox{ for all small }s)=1.
$$
The reader will have no trouble finding  explicit expressions for the transition
densities $p_t(x,y)$ and $q_t(x,y)$, thereby verifying that for $t>0$, $y>0$,
$$
\eqalign{
p_t(0+,y)=q_t(0-,y)&={1+e^{-2t}\over\sqrt{2\pi t^3}}e^{-y^2/2t}\cr&>
{1-e^{-2t}\over\sqrt{2\pi t^3}}e^{-y^2/2t}=p_t(0-,y)=q_t(0+,y),\cr}
$$
which is consistent with Theorem (1.16).
\medskip

This example is typical of what can go wrong when the hypothesis
[$P^{x,x}_{t_0}=Q^{x,x}_{t_0}$, $\forall\,x$] of Theorem (1.10) is weakened to
$P^{x_0,y_0}_{t_0}=Q^{x_0,y_0}_{t_0}$.
In general, under this latter condition, there is a set $N\in\EE$ that is both $X$-polar
and $Y$-polar and a set $\hat N\in\EE$ that is both $\hat X$-polar and $\hat
Y$-polar, such that the conclusions drawn in Theorem (1.10) remain true {\it
provided\/} one substitutes ``$x\in E\setminus N$'' for ``$x\in E$'' and ``$y\in
E\setminus
\hat N$'' for ``$y\in E$'' throughout. (Actually, the functions $\psi$ and $\hat\psi$ can
be defined so that (1.11) and (1.12) hold on all of $E$; these functions will be strictly
positive on $E$, but their finiteness can be guaranteed only off $N$ and $\hat N$,
respectively.) Since the proof of this assertions is quite close to  that of
Theorem (1.10), it is omitted.

After discussing bridge laws in section 2, we turn
to the proof of Theorem (1.10) in section 3. Theorem (1.16) is proved in section 4. 
\bigskip

\nin{\bf 2. Bridges}
\medskip

\nin The  discussion in this section is phrased in terms of  $X$, but applies equally
to $Y$. The process $X$ is as described in section 1.  All of the material in this
section, with the exception of Lemmas (2.8) and (2.9), is drawn from
\ref{\FPY}, to which we refer the reader for proofs and further discussion.

The following simple lemma shows that in constructing $P^{x,y}_t$ it matters not whether
we condition
$P^x$ on the event $\{X_t=x\}$ or on the event $\{X_{t-}=x\}$. 

\proclaim{(2.1) Lemma} $P^x(X_{t-}=X_t)=1$ for every $x\in E$ and every $t>0$.
\endproclaim

In what follows,  $\FF_{t-}$ denotes the $\sigma$-algebra generated by
$\{X_s, 0\le s<t\}$.

\proclaim{(2.2) Proposition} Given $(x,t,y)\in E\times ]0,\infty[\times E$  there is a
unique probability measure
$P^{x,y}_t$ on
$(\Omega,\FF_{t-})$ such that 
$$ 
\IPxyl(F)=P^x\left(F\cdot
{p_{t-s} (X_s,y)\over p_t(x,y)}\right)
\leqno(2.3) 
$$ 
for all positive $\FF_s$-measurable
functions
$F$ on $\Omega$, for all $0\le s<t$. 
Under $\IPxyl$ the coordinate process
$(X_s)_{0\le s< t}$ is a non-homogeneous strong Markov process with
transition densities 
$$
p^{(y,t)}(z,s;z',s')={p_{s'-s}(z,z')p_{t-s'}(z',y)\over
p_{t-s}(z,y)},\qquad 0<s<s'<t.
\leqno(2.4) 
$$ 
Moreover
$\IPxyl(X_0=x,X_{t-}=y)=1$. Finally, if $F\ge 0$ is
$\FF_{t-}$-measurable, and $g\ge0$ is a Borel function on $E$, then 
$$
P^x(F\cdot g(X_{t-}))=\int_E\IPxyl(F)\,g(y)\,p_t(x,y)\,m(dy).
\leqno(2.5) 
$$ 
Thus $(\IPxyl)_{y\in
E}$ is a regular version of the family of conditional probability distributions
$\{P^x(\cdot\,|X_{t-}=y)$, $y\in E\}$; equally so with $X_{t-}$ replaced by $X_t$,
because of Lemma (2.1).
\endproclaim

The following corollaries of Proposition (2.2) will be used in the sequel.

\proclaim{(2.6) Corollary} The $P^{x,y}_t$-law
of the time-reversed process $(X_{(t-s)-})_{0\le s< t}$ is
$\hat P^{y,x}_t$, the law of a $(y,t,x)$-bridge for the dual process
$\hat X$. \endproclaim

\proclaim{(2.7) Corollary}
For each $( \FF_{t+})$ stopping time $T$,
a
$\IPxyl$
regular conditional distribution for
$(X_{T + u } , 0 \le u < t - T ) $
given $\FF_{T+}$ on $\{T < t\}$ is provided by $P^{X_T ,y}_{t - T }$.
\endproclaim

Continuity properties are useful in trying to minimize the exceptional sets involved in
statements concerning bridge laws. The following simple result will be used in the proof
of (1.16).

\proclaim{(2.8) Lemma} Assume that  $x\mapsto p_t(x,y)$ is
 continuous  for each fixed pair $(t,y)\in]0,\infty[\times E$. Fix
$0<s<t$ and let $G$ be a bounded $\FF_{(t-s)-}$-measurable function on $\Omega$. Then
for each $y\in E$,
$$
x\mapsto P^{x,y}_{t}(G\comp\theta_s)
$$
is continuous on $E$.\endproclaim

\prf By Corollary (2.7),
$$
P^{x,y}_{t}(G\comp\theta_s)=\int_E {p_s(x,z)p_{t-s}(z,y)\over
p_{t}(x,y)}\,P^{z,y}_{t-s}(G)\,m^X(dz).
\leqno(2.9)
$$
The ratio on the
right side of (2.9) (call it $f_{x}(z)$) is a probability density with respect to
$m^X(dz)$, and the mapping $x\mapsto f_{x}(z)$ is continuous by hypothesis. It therefore
follows from Scheff\'e's Theorem \nref{\BIL}{p.\ 224} that
$x\mapsto f_{x}$ is a continuous mapping of $E$ into
$L^1(m^X)$. \qed
\bigskip

The {\it backward space-time process\/} associated with $X$ is the (Borel right) process
$$
\overline X_t(\omega,r):=(X_t(\omega),r-t),
$$
realized on the sample space
$\Omega\times{\bf R}$ equipped with the laws $P^x\otimes
\epsilon_r$.  A
(universally measurable) function $f:E\times{\bf R}\to[0,\infty]$ is $\overline
X$-excessive if and only if
$$
t\mapsto \int_E p_t(x,y) f(y,r-t)\, m^X(dy)
$$
is decreasing and right-continuous on $[0,\infty[$ for each $(x,r)\in E\times{\bf R}$.
For example, if $(y,s)\in E\times{\bf R}$ is fixed, then 
$(x,r)\mapsto 1_{]s,\infty[}(r)\,p_{r-s}(x,y)$ is ${\overline X}$-excessive.  A 
Borel function $f:E\times{\bf R}\to\overline{\bf R}$ is finely continuous with respect to
$\overline X$ if and only if $t\mapsto f(X_t,r-t)$ is right-continuous
$P^x\otimes\epsilon_r$-a.s.\ for every $(x,r)\in E\times{\bf R}$. Since
$\overline X$ is a right process \nref{\SHA}{\S 16}, $\overline X$-excessive functions are
finely continuous. Because of Hypotheses (1.4) and (1.6), the measure $m^X\otimes\leb$ on
$E\times{\bf R}$ is a reference measure for $\overline X$. Thus, if two finely continuous
functions of $\overline X$ agree $m^X\otimes\leb$-a.e., then they agree on all of
$E\times{\bf R}$. 

\proclaim{(2.10) Lemma} Fix $n\in{\bf N}$ and let $f_1, f_2,\ldots f_n$ be bounded
real-valued Borel functions on $E\times [0,\infty[$. Then for each $y\in E$, the function
$$
(x,t)\mapsto 1_{]0,\infty[}(t)\,P^{x,y}_t\left(\prod_{i=1}^n\int_0^t
f_i(X_s,t-s)\,ds\right)
\leqno(2.11)
$$
is finely continuous with respect to the backward space-time process
$(X_t,r-t)_{t\ge 0}$.\endproclaim

\prf Without loss of generality, we assume that $0<f_i\le 1$ for every $i$. 
The expression appearing in (2.11) can be written as the sum of $n!$ terms of the form
$$
1_{]0,\infty[}(t)P^{x,y}_t\int_0^t ds_1\int_{s_1}^t ds_2\cdots\int_{s_{n-1}}^t ds_n\,
\prod_{i=1}^n g_i(X_{s_i},t-s_i) ,
\leqno(2.12)
$$
where $(g_1,g_2,\ldots,g_n)$ is a permutation of $(f_1,f_2,\ldots,f_n)$. Let $h(x,t)$
denote the expression  in (2.12) multiplied by $p_t(x,y)$. Also, let 
$\tilde h(z,u):=p_u(z,y)\cdot P^{z,y}_u(J_u)$, where
$$
J_u:=\int_0^u du_2\int_{u_2}^u du_3\cdots\int_{u_{n-1}}^u du_n\,
\prod_{i=2}^n g_i(X_{u_i},u-u_i).
$$
For $t>0$, the Markov property (2.7) yields
$$
\eqalign{
h(x,t)&=p_t(x,y)\cdot P^{x,y}_t\int_0^tg_1(X_{s_1},t-s_1)
J_{t-{s_1}}\comp\theta_{s_1}\,ds_1\cr &=p_t(x,y)\cdot P^{x,y}_t\int_0^tg_1(X_{s_1},t-s_1)
P^{X(s_1),y}_{t-{s_1}}(J_{t-{s_1}})\,ds_1\cr
&=\int_E\int_0^t p_{s_1}(x,z) g_1(z,t-s_1) \tilde h(z,t-s_1)\,ds_1\,m^X(dz)\cr
&=\int_E\int_0^t p_{t-s}(x,z) g_1(z,s) \tilde h(z,s)\,ds\,m^X(dz).\cr}
\leqno(2.13)
$$
The final line in (2.13) exhibits $h$ as a positive linear combination of the space-time
excessive functions $(x,t)\mapsto 1_{]s,\infty[}(t)p_{t-s}(x,z)$, showing that $h$
is space-time excessive. Since $(x,t)\mapsto 1_{]0,\infty[}(t) p_t(x,y)$ is also
space-time excessive, the function appearing in (2.12) is finely continuous as
asserted.\qed
\bigskip

\nin{\bf 3. Proof of (1.10)}
\medskip

\nin For typographical convenience, throughout this section
 we assume (without loss of generality) that $t_0=2$, so  the basic hypothesis under
which we are working is that $Q^{x,x}_2=P^{x,x}_2$ for all $x\in E$.
\medskip

\nin{\it Proof of\/} (1.10)(a). Given $x\in E$ and $t\in]0,2[$, the mutual absolute
continuity of $P^x|_{\FF_t}$ and $Q^x|_{\FF_t}$ follows immediately from the hypothesis
$Q^{x,x}_{2}=P^{x,x}_{2}$ because of (2.3). Let us now show that if
$P^x|_{\FF_t}\sim Q^x|_{\FF_t}$ for all $x$, then $P^x|_{\FF_{2t}}\sim Q^x|_{\FF_{2t}}$
for all
$x$; an obvious induction will then complete the proof. By an application of the monotone
class theorem, given a bounded $\FF_{2t}$-measurable function
$F$, there is a bounded $\FF_t\otimes\FF_t$-measurable function $G$ such that
$F(\omega)=G(\omega,\theta_t\omega)$ for all $\omega\in\Omega$. Consequently,
$$
P^x(F)=\int_\Omega\int_\Omega
P^x(d\omega)\,P^{\omega(t)}(G(\omega,\cdot\,))
$$
and
$$
Q^x(F)=\int_\Omega\int_\Omega
Q^x(d\omega)\,Q^{\omega(t)}(G(\omega,\cdot\,))
$$
so the equivalence of $P^x$ and $Q^x$ on $\FF_{2t}$ follows from their equivalence
on $\FF_t$, as desired. The dual assertion can be proved in the same way once we
notice that $\hat Q^{x,x}_2=\hat P^{x,x}_2$ for all $x\in E$, because of Corollary
(2.6).\qed
\medskip

An important consequence of the equivalence just proved is that $X$ and $Y$ have the same
fine topologies, as do their space-time processes. Of course, the same can be said of
$\hat X$ and $\hat Y$.
\medskip

\nin{\it Proof of\/} (1.10)(b).   The argument is broken into several steps.

{\it Step 1:\/} $m^X\sim m^Y$. Indeed, because the transition densities are strictly
positive and finite by hypothesis, $m^X$ is equivalent to the
$P^{x,x}_2$-distribution of $X_{1}$, while $m^Y$ is equivalent to the
$Q^{x,x}_2$-distribution of $Y_{1}$ (for any fixed $x\in E$).

{\it Step 2.\/} {\sl For each $(x,t)\in E\times]0,2[$, $Q^{x,y}_t=P^{x,y}_t$ for
$m^X$-a.e.\ $y\in E$.\/} Fix $(x,t)\in E\times]0,2[$. Then by (2.6) and (2.7), the
$P^{x,x}_2$-conditional distribution of
$(X_s)_{0\le s<t}$, given $X_{t-}=y$, is $P^{x,y}_t$ (for $m^X$-a.e.\ $y\in E$).
Similarly, the $Q^{x,x}_2$-conditional
distribution of
$(Y_s)_{0\le s<t}$, given $Y_{t-}=y$, is $Q^{x,y}_t$ (for $m^Y$-a.e.\ $y\in E$).
The assertion therefore follows from the basic hypothesis ($Q^{x,x}_2=P^{x,x}_2$, $\forall
x$) because of Step 1.

{\it Step 3.\/} {\sl There exists $b\in
E$ such that
$Q^{x,b}_t=P^{x,b}_t$ for all $x\in E$ and all $t\in]0,2[$.\/}  
By Step 2 and Fubini's theorem there exists $b\in E$ such that $P^{x,b}_t=Q^{x,b}_t$
for $m^X\otimes\leb$-a.e.\ $(x,t)\in E\times ]0,2[$. Let
$\II$ denote the class of processes
$I$ of the form
$$
I_t:=\prod_{i=1}^n\int_0^t f_i(X_s,t-s)\,ds,\qquad t\ge 0,
$$
where $n\in{\bf N}$ and each $f_i$ is a bounded real-valued Borel function on
$E\times[0,\infty[$. It is easy to see that for each fixed $t>0$, the 
family
$\{I_t:I\in\II\}$ is measure-determining on  $(\Omega,\FF_{t-})$. Therefore, it suffices
to show that 
$$
P^{x,b}_t(I_t)=Q^{x,b}_t(I_t)
\leqno(3.1)
$$ 
for all $x\in E$,
$t\in]0,2[$, and $I\in\II$. But by Lemma (2.10) and the remark made following the proof
of (1.10)(a), the two sides of (3.1) are finely-continuous (with respect to the space-time
processes
$(X_t,r-t)_{t\ge 0}$ and
$(Y_t,r-t)_{t\ge 0}$) on all of
$E\times]0,\infty[$, as functions of
$(x,t)$.  By the choice of $b$ these  functions agree
$m^X\otimes\leb$-a.e.\ on the (space-time) finely open set $E\times]0,2[$; consequently, 
they agree everywhere on $E\times]0,2[$, because $m^X\otimes\leb$ is a reference measure
for the space-time processes.

{\it Step 4.\/} In view of Step 3 there exists $b\in E$ such that $P^{x,b}_1=Q^{x,b}_1$
for all $x\in E$. This $b$ will remain fixed in the following discussion. Recall from
(1.10)(a) that the laws $P^x$ and $Q^x$ are (locally) mutually absolutely continuous
for each $x\in E$. Let
$Z_t$ denote the Radon-Nikodym derivative
$dP^x|_{\FF_{t+}}/dQ^x|_{\FF_{t+}}$. Then $Z$ is a strictly positive
right-continuous martingale and a multiplicative functional of $X$; see, for example,
\nref{\KUN}{Thm.\ 5.1}. The term {\it multiplicative\/} refers to the identity
$$
Z_{t+s}=Z_t\cdot Z_s\comp\theta_t,\qquad P^x\hbox{-a.s.},\,\forall\,x\in E,
\forall\,s,t\ge 0.
$$
Using (2.3) we see that for any $x\in E$,
$$
\eqalign{
P^{x,b}_1(F)
&=Q^{x,b}_1(F)=Q^x\left(F\,{q_{1-s}(X_s,b)\over q_1(x,b)}\right)\cr
&=P^x\left(F\cdot Z_s\,{q_{1-s}(X_s,b)\over q_1(x,b)}\right)\cr
&=P^{x,b}_1\left(F\cdot Z_s\,{q_{1-s}(X_s,b)\over q_1(x,b)}{p_{1}(x,b)\over
p_{1-s}(X_s,b)}\right)\cr }
$$
for any $F\in\FF_{s+}$, provided $0<s<1$. Since $Z_s$ is $\FF_{s+}$ measurable, it follows
that
$$
Z_s={p_{1-s}(X_s,b)\over q_{1-s}(X_s,b)}{q_{1}(x,b)\over p_1(x,b)}\qquad
P^{x,b}_1\hbox{-a.s.}
\leqno(3.2)
$$
for all $x\in E$ and $0<s<1$. 
Since $P^x$ and $P^{x,b}_1$ are
equivalent on
$\FF_{s+}$ for $0<s<1$, we see that
$$
Z_s=\varphi_s(X_0,X_s)\qquad 
P^x\hbox{-a.s.}, \forall\,s \in]0,1[, \forall x\in E,
$$
where 
$$
\varphi_s(x,z):={\psi_s(z)\over\psi_0(x)}
$$
and
$$
\psi_s(z):={p_{1-s}(z,b)\over q_{1-s}(z,b)}.
$$
The function $(z,s)\mapsto 1_{[0,1[}(s)p_{1-s}(z,b)$ is an
excessive function of the {\it forward\/} space-time process $(X_t,t+r)_{t\ge 0}$
restricted to $E\times[0,\infty[$; it is therefore space-time finely continuous on
$E\times[0,1[$. In the same way
$(z,s)\mapsto q_{1-s}(z,b)$ is finely continuous on $E\times[0,1[$ with
respect to the space-time process
$(Y_t,t+r)_{t\ge 0}$. But the fine topology of the latter process is the same as that of
$(X,r+t)_{t\ge 0}$ because of the mutual absolute continuity ($P^x|_{\FF_t}\sim
Q^x|_{\FF_t}$, $\forall (x,t)$)  already established. It follows that $(z,s)\mapsto
\psi_s(z)$ is space-time finely continuous on
$E\times [0,1[$.
Now from  the multiplicativity  of $Z$ and
the strict positivity of the transition densities of $X$ we deduce that for all $x\in E$
and all $t,s>0$ such that $t+s<1$, there is an $m^X\otimes m^X$-null set $N(x,t,s)\subset
E\times E$ such that
$$
\varphi_{t+s}(x,y)=\varphi_t(x,z)\cdot\varphi_s(z,y)
\leqno(3.3)
$$
provided $(y,z)\notin N(x,t,s)$. By the preceding discussion, the two sides of (3.3) are
space-time finely continuous as functions of $(y,s)$. Moreover, 
$m^X\otimes\leb$  is a reference measure for $(X_t,r+t)$; thus, two space-time
finely continuous functions equal $m^X\otimes\leb$-a.e.\ must be identical. From this
observation and Fubini's theorem it follows that given
$(x,t)\in E\times ]0,1[$ there is an $m^X$-null set $N(x,t)$ such that (3.3) holds for
all $(y,s)\in E\times[0,1-t[$ and all $z\notin N(x,t)$. Taking $s=0$ we find that
$$
{\psi_0(y)\over\psi_t(y)}={\psi_0(z)\over\psi_t(z)}
\leqno(3.4)
$$
for all $y\in E$, $0<t<1$, and $z\notin N(x,t)$. Thus, defining
$\lambda_t:=-\log[\psi_t(b)/\psi_0(b)]$ and $\psi:=\psi_0$, we have, for each $x\in E$,
$$
Z_t=e^{-\lambda_t}\,{\psi(X_t)\over\psi(X_0)},\qquad P^x\hbox{-a.s.},
\leqno(3.5)
$$
for all $t\in]0,1[$, since $P^x(X_t\in N)=0$ for any $m^X$-null set $N$. The
multiplicativity of
$Z$ implies first that $\lambda_t=\lambda t$ for some real constant $\lambda$, and then
that (3.5) holds for all $t> 0$. This yields (1.13), from
which (1.11)   follows immediately because $Z$ is a
$P^x$-martingale.

The dual assertions (1.12) and (1.14) are proved in the same way, and the fact that
$\psi$ and $\hat\psi$ correspond to the same ``eigenvalue'' $\lambda$ follows easily
from (1.5). 

Turning to the final assertion, let $\rho$ denote a strictly positive and finite
version of the Radon-Nikodym derivative $dm^Y/dm^X$---the equivalence of $m^X$ and
$m^Y$ follows immediately from (1.13). Using (1.5) (for
$X$ and for $Y$) one can check that $P_t(\rho/\psi\hat\psi)=\rho/\psi\hat\psi$ and
$P_t(\psi\hat\psi/\rho)=\psi\hat\psi/\rho$, $m^X$-a.e. 
Consequently,
$$
1=P_t1=P_t\Big((\rho/\psi\hat\psi)^{1/2}(\psi\hat\psi/\rho)^{1/2}\Big)
\le \Big( P_t(\rho/\psi\hat\psi)\,P_t(\psi\hat\psi/\rho)\Big)^{1/2}=1,
$$
which forces  $\rho=\psi\hat\psi$, $m^X$-a.e, as claimed.
\qed
\medskip

\nin{\it Proof of\/} (1.10)(c). Formula (1.13) implies that for each $x\in E$ and
$t>0$,
$$
q_t(x,y)=e^{-\lambda t}{1\over\psi(x)\hat\psi(y)}p_t(x,y),\qquad m^X\hbox{-a.e.\
}y\in E,
\leqno(3.6)
$$
because $\psi\hat\psi=dm^Y/dm^X$.
For fixed $x$ the two sides of (3.6) are finely continuous (as functions of
$(y,t)\in E\times ]0,\infty[$) with respect to the  backward space-time
process
$(\hat X_t,r-t)_{t\ge 0}$. (As before, the equivalence of laws established in (1.10)(a)
implies that $(\hat X_t,r-t)$ and $(\hat Y_t,r-t)$ have the same fine topologies.) Since
$m^X\otimes\leb$ is a reference measure for this space-time process, the equality in
(3.6) holds for {\it all\/} $(y,t)\in E\times]0,\infty[$. The asserted equality of bridges
now follows from (1.13) and (2.3).\qed
\bigskip

\nin{\bf 4. Proof of  (1.16)}
\medskip

\nin We first show that
$P^{x,y_0}_{t_1}=Q^{x,y_0}_{t_1}$ for all $x\in E$, where $t_1:=t_0/2$. To this end fix
$x\in E$, let
$d$ be a metric on
$E$ compatible with its topology, and let
$B(\delta)$ denote the $d$-ball of radius $\delta$ centered at $x$. Let $F$ be a bounded
$\FF_{t_1-}$-measurable function of the form $G\comp\theta_s$, where $0<s<t_1$ and
$G\in\FF_{(t_1-s)-}$. By Corollary (2.7),
$$
\eqalign{
P^{x_0,y_0}_{t_0}(F\comp\theta_{t_1}\,&|X_{t_1}\in B(\delta))\cr
&=\int_{{B(\delta)}}\,P^{z,y_0}_{t_1}(F)\,P^{x_0,y_0}_{t_0}(X_{t_1}\in
dz)/P^{x_0,y_0}_{t_0}(X_{t_1}\in B(\delta)).\cr}
\leqno(4.1)
$$
(Notice that $P^{x_0,y_0}_{t_0}(X_{t_1}\in B(\delta))>0$ because of the strict
positivity of the transition density function of $X$.)  By Lemma (2.8), the mapping
$z\mapsto P^{z,y_0}_{t_1}(F)$ is continuous. Since  the probability measure
$$
dz\mapsto 1_{B(\delta)}(z)\,P^{x_0,y_0}_{t_0}(X_{t_1}\in
dz)/P^{x_0,y_0}_{t_0}(X_{t_1}\in B(\delta))
$$
converges weakly to the unit mass at $x$ as $\delta\to 0$, it follows from (3.7) that
$$
P^{x,y_0}_{t_1}(F)=\lim_{\delta\to 0}P^{x_0,y_0}_{t_0}(F\comp\theta_{t_1}\,|\,X_{t_1}\in
B(\delta)).
\leqno(4.2)
$$
By hypothesis, the right side of (4.2) is unchanged if  $P^{x_0,y_0}_{t_0}$ is
replaced by $Q^{x_0,y_0}_{t_0}$; the same is therefore true of the left side, so 
$P^{x,y_0}_{t_1}(F)=Q^{x,y_0}_{t_1}(F)$. The monotone class theorem clinches the
matter.

The arguments used in the proof of Theorem (1.10) (especially Step 4 of the proof of
(1.10)(b)) can now be used to finish the proof. The dual assertion follows in the same
way. \qed

\vfill\eject

\frenchspacing
\baselineskip=13pt
\centerline{\bf References}
\bigskip

\itemitem{[\BL]}
Benjamini, I.\ and Lee S.: Conditioned diffusions which are Brownian motions. {\it J.\
Theor. Probab.\/} {\bf 10} (1997) 733--736.
\smallskip

\itemitem{[\BIL]}
Billingsley, P.: {\it Convergence of Probability Measures.} Wiley, New York, 1968.
\smallskip

\itemitem{[\DYN]}
E.B.\ Dynkin: Minimal excessive measures and functions. {\it Trans. Amer. Math. Soc.}
{\bf 258} (1980) 217--244. 
\smallskip

\itemitem{[\FPY]}
Fitzsimmons, P.J., Pitman, J.\ and Yor, M.: Markovian bridges: construction, Palm
interpretation, and splicing. In {\it Seminar on Stochastic Processes, 1992,\/} pp.\
101--133. Birkh\"auser, Boston, 1993.
\smallskip

\itemitem{[\FOL]}
F\"ollmer, H.: Martin boundaries on Wiener space. In {\it Diffusion processes and related
problems in analysis, Vol. I\/}, pp.\ 3--16.  BirkhŠuser Boston, Boston,  1990.
\smallskip

\itemitem{[\GET]}
Getoor, R.K.: {\it Markov processes: Ray processes and right processes.\/} Lecture Notes
in Math. {\bf 440}. Springer-Verlag, Berlin-New York, 1975.
\smallskip

\itemitem{[\GSexc]} 
Getoor, R.K.\ and Sharpe, M.J.: Excursions of dual processes. {\it Adv. Math.} {\bf 45}
(1982) 259--309.
\smallskip

\itemitem{[\GSnswd]} 
Getoor, R.K.\ and Sharpe, M.J.: Naturality, standardness, and weak duality for Markov
processes. {\it Z.\ Wahrscheinlichkeitstheorie verw.\ Gebiete\/} {\bf 67} (1984) 1--62.
\smallskip

\itemitem{[\KW]}
Kunita, H.\ and  Watanabe, T.: Markov processes and Martin boundaries, I.
{\it Illinois J. Math.\/} {\bf 9} (1965) 485--526.
\smallskip

\itemitem{[\KUN]}
Kunita, H.:    Absolute continuity of Markov processes. In {\it S\'eminaire de
Probabilit\'es X.\/} pp.\ 44--77. Lecture Notes in
Math. {\bf 511}, Springer, Berlin, 1976. 
\smallskip

\itemitem{[\SHA]}
Sharpe, M.J.: {\it General Theory of Markov Processes.\/} Academic Press, San Diego, 1988.
\smallskip

\itemitem{[\WIT]}
R.\ Wittmann: Natural densities of Markov transition probabilities. 
{\it Prob. Theor. Rel. Fields} {\bf 73} (1986) 1--10. 
\smallskip

\itemitem{[\YAN]}
J.-A.\ Yan: A formula for densities of transition functions. 
In {\it S\'eminaire de Probabilit\'es XXII,\/} pp.\ 92--100.  Lecture Notes in Math. {\bf
1321}, Springer, Berlin, 1988.
\smallskip

\end